\documentclass[a4paper]{article}

\author{Maria Alberich-Carrami\~nana --
        Joaquim Ro\'e
   \thanks {Both authors were partially supported by
CAICYT PB98-1185, Generalitat de Catalunya 98-SGR-00024
and ``EAGER'', European Union contract HPRN-CT-2000-00099.} \\
  \small{Departament d'\`Algebra i Geometria,
           Universitat de Barcelona,} \\
  \small{Gran Via, 585, E-08007, Barcelona. }\\
  \small{e-mail: alberich@mat.ub.es, jroevell@mat.ub.es} \\
           }

\usepackage{amsmath}
\usepackage{amsthm}
\usepackage{amssymb}
\usepackage{xypic}
\usepackage{psfig}

\theoremstyle{definition}
\newtheorem{Def}{Definition}[section]
\theoremstyle{plain}
\newtheorem{Lem}[Def]{Lemma}
\newtheorem{Cor}[Def]{Corollary}
\newtheorem{Pro}[Def]{Proposition}
\newtheorem{Teo}[Def]{Theorem}
\theoremstyle{remark}
\newtheorem{Rem}[Def]{Remark}
\newtheorem{Exa}{Example}

\newcommand{\Hilb}{\operatorname{Hilb}}
\newcommand{\Spec}{\operatorname{Spec}}

\newcommand{\Cl}{{\mathit{Cl}}}

\newcommand{\D}{{\bf D}}

\newcommand{\C}{\mathbb{C}}
\newcommand{\Z}{\mathbb{Z}}

\newcommand{\nnu}{{\boldsymbol \nu}}

\begin{document}

\title{Enriques diagrams and adjacency of planar curve singularities}
\maketitle

\begin{abstract}
We study adjacency of equisingularity types of planar
curve singularities in terms of their Enriques diagrams.
The goal is, given two equisingularity types,
to determine whether one of them is adjacent to the other.
For linear adjacency a complete answer is
obtained, whereas for arbitrary (analytic) adjacency
a necessary condition and a sufficient condition
are proved.
We also show an example
of singular curve of type $\mathcal{D}'$
that can be deformed to a curve of type
$\mathcal{D}$ without $\mathcal{D}'$
being adjacent to $\mathcal{D}$.
This suggests that
\emph{analytical} rather than topological
equivalence should be considered when
studying adjacency of singularity types.
\end{abstract}

\section*{Introduction}
A class of reduced (germs of) planar curve singularities
$\mathcal{D}'$ is said to be \emph{adjacent}
to the class $\mathcal{D}$ when every member of the
class $\mathcal{D}'$ can be deformed
into a member of the class $\mathcal{D}$ by an arbitrarily
small deformation. If this can be done with a linear deformation,
then we say that $\mathcal{D}'$ is \emph{linearly adjacent}
to $\mathcal{D}$.
It is well known that equisingularity and topological equivalence of
reduced germs of curves on smooth surfaces are equivalent, and that
analytical equivalence implies topological equivalence (see for
instance \cite{Zar32}, \cite{Zar65I} or \cite{Cas00}).
We shall focus on the
equisingularity (or topological equivalence) classes,
and we will call them simply \emph{types}. The Enriques diagrams
introduced by Enriques in \cite[IV.I]{EC15} represent the types:
two reduced curves are equisingular at $O$ if and only if
their associated Enriques diagrams are isomorphic
(see \cite[3.9]{Cas00}).

In \cite{Arn76} Arnold classified critical points of
functions with modality at most two; this implies the
classification of types (of planar curve singularities)
with multiplicity at most four. He also described some
adjacencies between them, introducing the so-called
\emph{series} of types $A$, $D$, $E$ and $J$. The
construction of series was generalized by Siersma in
\cite{Sie77} using a kind of Enriques diagrams;
in particular, he classified types of multiplicity
at most five. As we shall see below, all adjacencies
within one series are linear.

Apart from series, only some particular cases of adjacency
are known, obtained using explicit deformations.
On the other hand, the semicontinuity of some numerical
invariants such as the genus discrepancy $\delta$ or
the Milnor number $\mu$ provide necessary conditions
for adjacency, but these invariants are far from
determining the type of a singularity, and so it
is not to be expected that they give a complete answer
to the adjacency question.
Here, instead of numerical invariants, the Enriques
diagram (which does determine the type) is used,
providing a necessary condition and a sufficient
condition for adjacency.
In the case of \emph{linear} adjacency
a complete answer follows, namely, we determine all
linear adjacencies in terms of Enriques diagrams.
Non-linear adjacencies are a much subtler subject, as
shown by the fact that the types do not form a stratification
of $\C[[x,y]]$ (see example \ref{contrex}); this
suggests that a complete understanding
of analytic adjacencies can only be achieved by
considering analytic moduli of singularities,
rather than equisingularity classes alone.
%\enlargethispage{\baselineskip}

We present a definition of Enriques diagrams in a purely combinatorial
way, that was used by Kleiman and Piene (\cite{KP99})
to list all equisingularity types with codimension up to 8, which
is needed for the enumeration of 8--nodal curves
(see also \cite{GSG92}).

A \emph{tree} is a finite directed graph, without
loops; it has a single initial vertex, or \emph{root}, and every other
vertex has a unique immediate predecessor.
If $p$ is the immediate predecessor of the vertex
$q$, we say that $q$ is a successor of $p$. If $p$
has no successors then it is an extremal vertex.
An \emph{Enriques diagram} is a tree
with a binary relation between vertices, called \emph{proximity},
which satisfies:
\begin{enumerate}
\item The root is proximate to no vertex.
\item Every vertex that is not the root is proximate to
   its immediate predecessor.
\item No vertex is proximate to more than two vertices.
\item If a vertex $q$ is proximate to two vertices
   then one of them is the immediate predecessor of $q$,
   and it is proximate to the other.
\item Given two vertices $p$, $q$ with $q$ proximate to $p$,
   there is at most one vertex proximate to both of them.
\end{enumerate}
The vertices which are proximate to two points are called
\emph{satellite}, the other vertices are called \emph{free}.
We usually denote the set of vertices of an Enriques diagram
$\D$ with the same letter $\D$.

To show graphically the proximity relation, Enriques diagrams are
drawn according to the following rules:
\begin{enumerate}
\item If $q$ is a free successor of $p$
   then the edge going from $p$ to $q$ is smooth and curved and,
   if $p$ is not the root, it has at $p$ the same tangent
   as the edge joining $p$ to its predecessor.
\item The sequence of edges connecting a maximal
   succession of vertices proximate to the same vertex $p$
   are shaped into a line segment, orthogonal to the edge joining $p$
   to the first vertex of the sequence.
\end{enumerate}
An \emph{isomorphism} of Enriques diagrams is a bijection $i$
between the sets of vertices of the two diagrams so that
$q$ is proximate to $p$ if and only if $i(q)$ is proximate to
$i(p)$; two Enriques diagrams are \emph{isomorphic} if there
is an isomorphism between them.

\begin{Exa}
\label{exa1}
Figure \ref{dia1} shows an Enriques diagram with
nine vertices. $p_1$ is the root of the
diagram, $p_4, p_5$ are satellites proximate to $p_2$,
$p_6$ is a satellite proximate to $p_3$ and the
remaining vertices are free.
\begin{figure}
  \begin{center}
    \mbox{\psfig{file=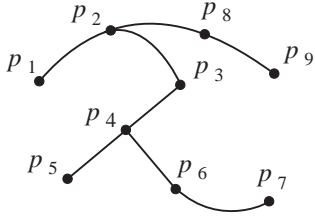}}
%    \mbox{\includegraphics{diagrama.eps}}
    \caption{Enriques diagram of example \ref{exa1}. \label{dia1}}
  \end{center}
\end{figure}
\end{Exa}

A \emph{subdiagram} of an Enriques diagram $\D$ is a subtree
$\D_0 \subset \D$
together with the induced proximity relation, such that the
predecessors of every vertex $q \in \D_0$ belong to $\D_0$.
An \emph{admissible ordering} for an Enriques diagram
$\D$ is a total ordering $\preceq$ for its set of vertices
refining the natural ordering of $\D$, i.e., such that for
every vertex $p$, and every successor $q$ of $p$, $p \preceq q$.

Given an Enriques diagram $\D$ of $n$ vertices with an admissible
ordering $\preceq$, let $p_1$, $p_2$, \dots, $p_n$ denote its vertices,
numbered according to $\preceq$.
The \emph{proximity matrix} of $\D$ is a square matrix $P=(p_{i,j})$
of order $n$, with
$$
p_{i,j}=
\begin{cases}
  1 & \text{if }i=j,\\
 -1 & \text{if $p_i$ is proximate to $p_j$},\\
  0 & \text{otherwise.}
\end{cases}
$$

A \emph{system of multiplicities} for (the vertices of) an Enriques diagram
$\D$ is any map $\nu:\D \rightarrow \Z$. We will usually write
$\nu_p=\nu(p)$. A pair $(\D,\nu)$, where $\D$ is an Enriques
diagram and $\nu$ a system of multiplicities for it, will be called a
\emph{weighted Enriques diagram}. The \emph{degree} of a
weighted Enriques diagram is
$\deg (\D,\nu)= \sum_{p \in \D} \nu_p(\nu_p+1)/2$.
A \emph{consistent Enriques diagram}
is a weighted Enriques diagram such that, for all $p\in \D$,
$$ \nu_p \geq \sum_{q \text{ prox. to } p}  \nu_q.$$
Note that if $(\D,\preceq)$ is an Enriques diagram of $n$ vertices
with an admissible ordering, then a system of multiplicities for $\D$
may be identified with a vector $\nnu=(\nu_1, \nu_2, \dots, \nu_n)\in
\Z^n$, taking $\nu_i=\nu_{p_i}$, $i=1, \dots, n$; we shall
use the notation $(\D,\preceq, \nnu)$ for a weighted ordered
Enriques diagram, where $\nnu \in \Z^n$.

\begin{figure}
  \begin{center}
    \mbox{\psfig{file=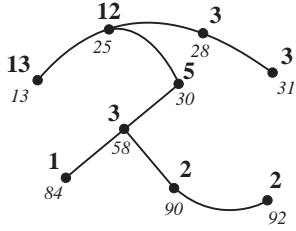}}
%    \mbox{\includegraphics{multval.eps}}
    \caption{A consistent system of multiplicities (in bold shape)
      for the Enriques diagram of example \ref{exa1},
      and its associated system of values (in italics). \label{dia2}}
  \end{center}
\end{figure}

To every system of multiplicities $\nu$ for a diagram $\D$ we
associate a \emph{system of values}, which is another map
$v:\D \rightarrow \Z$, defined recursively as
$$
v_p=
\begin{cases}
\nu_p & \text{if $p$ is the root,} \\
\nu_p + \sum_{p \text{ prox. to } q} v_q & \text{otherwise.}
\end{cases}
$$
Observe that any map
$v:\D \rightarrow \Z$ is the system of values associated the system of
multiplicities $\nu:\D \rightarrow \Z$ defined recursively as
$$
\nu_p=
\begin{cases}
v_p & \text{if $p$ is the root,} \\
v_p - \sum_{p \text{ prox. to } q} v_q & \text{otherwise.}
\end{cases}
$$
Hence giving a system of multiplicities for an Enriques diagram is
equivalent to give a system of values.
Figure \ref{dia2} shows the system of values
associated to a consistent system of multiplicities.

%\end{Def}

The relationship between the combinatorial properties of Enriques
diagrams
and the topology of planar curve singularities is explained next.
%In order to keep the context as in Arnold's paper, we
Assume that $O$ is a smooth point on a complex surface $S$,
whose local ring is isomorphic to $\C[[x,y]]$, and
let $f \in \C[[x,y]]$ be the equation of a (germ of)
curve with an isolated singularity at $O$.

Let $K$ be a finite set of points equal or infinitely near to
the smooth point $O$, such that for each $p\in K$, $K$ contains
all points to which $p$ is infinitely near. Such a set is
called a \emph{cluster} of points infinitely near to $O$.
A point $p\in K$ is said to be \emph{proximate} to another
$q \in K$ if it is infinitely near to $q$ and lies on the \emph{strict}
transform of the exceptional divisor of blowing up $q$. Thus
being proximate to is a binary relation between points of a cluster
which satisfies the same conditions as proximity between vertices of an
Enriques diagram.
Therefore for every cluster
there is an associated Enriques diagram which encodes all the
information on proximities between points of $K$.

The \emph{value} of a germ of curve at a point
$p$ of a cluster $K$ is the multiplicity at $p$ of the
pullback of the germ of curve in the blown up surface containing
$p$; in the case $p=O$ the value is just the multiplicity
of the curve at $O$.
Given a cluster $K$ and a system of values $v: K \rightarrow \Z$
(associated to the system of multiplicities $\nu $)
there is a complete ideal $H_{K,\nu} \subset \C[[x,y]]$ containing
all equations of the germs of curve which have at every point $p\in K$
\emph{value} at least $v_p$ (see \cite[4.4.4]{Cas00}).

%\enlargethispage{6mm}

\begin{Exa}
If $C \subset S$ is a reduced curve going through $O$, then the
set of singular points of $C$ equal or infinitely near to $O$
is a cluster $K$. The Enriques diagram of $K$, weighted
with the multiplicities of $C$ at the points of $K$, is a consistent
Enriques diagram, which we call the Enriques diagram associated
to the singularity of $C$. Such a diagram
has no extremal free vertices of multiplicity
$\nu_p \le 1$ (because $K$ consists only of singular points of $C$,
which either have multiplicity bigger than 1, are satellites
or precede some satellite point on $C$).
Conversely, if $\D$ is a consistent
Enriques diagram with no extremal free vertices of multiplicity
$\nu_p \le 1$ then there are germs of curve at $O$ whose cluster
of singular
points has Enriques diagram isomorphic to $\D$ (see \cite{Cas00}).
\end{Exa}

It is well known that two reduced curves are equisingular at $O$
if and only if
their associated Enriques diagrams
are isomorphic (see \cite[3.8]{Cas00}, for instance).
Therefore the equisingularity types (types for short) of reduced
germs of curves on smooth surfaces are identified
with the isomorphism classes of consistent
Enriques diagrams with no extremal free vertices of multiplicity
$\nu_p \le 1$.

\section{Linear adjacency}

Let $I \subset \C[[x,y]]$ be an ideal.
According to \cite[7.2.13]{Cas00} general members (by the Zariski
topology of the coefficients of
the series) of $I$ define
equisingular germs.

\begin{Lem}
$(\D',\nu')$ is linearly adjacent
to $(\D,\nu)$ if and only if for every $f \in \C[[x,y]]$
defining a reduced germ of curve of type $(\D',\nu')$,
there exists an ideal $I \subset \C[[x,y]]$ with
$f \in I$ and whose general member defines a reduced germ of type
$(\D,\nu)$.
\end{Lem}
\begin{proof}
The if part of the claim is evident. To see the only if
part, assume that $f$ defines a reduced germ of type $(\D',\nu')$
that can be deformed to a reduced germ of type $(\D,\nu)$
by a linear deformation $f+tg$, $g \in \C[[x,y]]$.
This means that general members of the pencil $f+tg$
define reduced germs of type $(\D,\nu)$.
Hence
general members of the ideal $I=(f,g)$
define germs of type $(\D,\nu)$ as well
(\cite[7.2]{Cas00}).
\end{proof}

\begin{Pro}
\label{cideal}
Let $(\D, \mu)$ and $(\D',\mu')$ be
weighted Enriques diagrams, with $(\D',\mu')$
consistent. The following are
equivalent:
\begin{enumerate}
\item \label{ex}There are two clusters,
$K$ and $K'$, whose Enriques diagrams are
$\D$ and $\D'$ respectively, such that
$H_{K',\mu'} \subseteq H_{K,\mu}$.
\item \label{bs}For every cluster $K$ with
Enriques diagram $\D$, there is a cluster
$K'$ with Enriques diagram $\D'$ such that
$H_{K',\mu'} \subseteq H_{K,\mu}$.
\item \label{sb}For every cluster $K'$ with
Enriques diagram $\D'$, there is a cluster
$K$ with Enriques diagram $\D$ such that
$H_{K',\mu'} \subseteq H_{K,\mu}$.
\item \label{comb}There exist
isomorphic subdiagrams $\D_0 \subset \D$,
$\D_0' \subset \D'$ and an isomorphism
$$ i: \D_0 \longrightarrow \D_0'$$
such that the system of multiplicities
$\nu$ for $\D$ defined as
$$
\nu(p)=
\begin{cases}
\mu'(i(p)) & \text{ if } p \in \D_0 , \\
0 & \text{ otherwise}
\end{cases}
$$
has the property that the values $v$ and $v'$ associated
to the multiplicities $\mu$ and $\nu$ respectively
satisfy $v(p) \le v'(p) \ \forall p \in \D$.
\end{enumerate}
\end{Pro}

\begin{proof}
Clearly both \ref{sb} and \ref{bs} imply
\ref{ex}. We shall prove that \ref{ex}
implies \ref{comb} and that \ref{comb}
implies both \ref{bs} and \ref{sb}.

Let us first prove that \ref{ex} implies \ref{comb}.
So assume there are two clusters,
$K$ and $K'$, whose Enriques diagrams are
$\D$ and $\D'$ respectively, such that
$H_{K',\mu'} \subseteq H_{K,\mu}$.
The points common to $K$ and $K'$ clearly
form a cluster, which we call $K_0$.
The vertices in $\D$ and $\D'$
corresponding to points in $K_0$ form
subdiagrams $\D_0$ and $\D_0'$, and
the coincidence of points in $K_0$
determines an isomorphism
$ i: \D_0 \longrightarrow \D_0'$.
It only remains to be seen that
the values $v$ and $v'$ associated
to the multiplicities $\mu$ and $\nu$ respectively
(with $\nu$ as in the claim)
satisfy $v(p) \le v'(p) \ \forall p \in \D$.
Now choose a germ $f \in H_{K',\mu'}$ having
multiplicity exactly $\mu_p'$ at each point
$p \in K'$ (such an $f$ exists because
$(\D',\mu')$ is consistent, see \cite[4.2.7]{Cas00}).
This implies that $f$ has value exactly
$v'(p)$ at each point $p\in K$.
Then $f \in H_{K,\mu}$ because
$H_{K',\mu'} \subseteq H_{K,\mu}$,
and the claim follows by the
definition of $H_{K,\mu}$.

Let us now prove that \ref{comb} implies \ref{sb}.
Assume that \ref{comb} holds, and let $K'$
be a cluster whose Enriques diagram is $\D'$. We
must prove the existence of a cluster $K$ with
Enriques diagram $\D$ such that
$H_{K',\mu'} \subseteq H_{K,\mu}$.
Let $K_0$ be the cluster
formed by the points corresponding to vertices in $\D_0'$.
Add to $K_0$ the points necessary to get a cluster $K$ with
Enriques diagram $\D$. Because of the
hypothesis on the values $v$ and $v'$
and the characterization of $H_{K,\mu}$ (see for
instance \cite[4.5.4]{Cas00}), $H_{K',\mu'} \subseteq H_{K,\mu}$.

In the same way it is proved that
\ref{comb} implies \ref{bs}.
\end{proof}

If the conditions of proposition \ref{cideal}
are satisfied, we shall write
$(\D', \mu') \ge(\D, \mu)$.
Now we can prove our main result.
The interest of proposition \ref{cideal}
and theorem \ref{linesp}
lies on the fact that condition \ref{comb} of \ref{cideal}
can be checked directly on the Enriques diagrams,
using their combinatorial properties, thus giving a practical
means to decide whether a type is
or is not linearly adjacent to another.

\begin{Teo}
\label{linesp}
Let $(\D,\mu)$, $(\tilde \D,\tilde \mu)$ be types.
$(\tilde \D,\tilde \mu)$ is linearly adjacent to $(\D,\mu)$ if
and only if there exists a weighted consistent Enriques diagram
$(\D',\mu')$, differing from $(\tilde \D,\tilde \mu)$ at most
in some free vertices of multiplicity one,
satisfying $(\D',\mu') \ge (\D,\mu)$.
\end{Teo}

\begin{proof}
To prove the if part,
given a reduced germ $f \in \C[[x,y]]$
defining a curve singularity of type $(\tilde \D, \tilde \mu)$
we have to show the existence of
an ideal $I \subset \C[[x,y]]$ containing $f$ and whose general
member defines a reduced germ of type $(\D, \mu)$,
provided that $(\D',\mu') \ge (\D,\mu)$.
Let $\tilde K$ be the cluster of singular points of $f$ (whose
Enriques diagram is $\tilde \D$). For each vertex $p$
of $\D'$ not in $\tilde \D$, whose predecessor is denoted by $q$,
choose a point on $f=0$ on the first neighbourhood
of the point corresponding to the vertex $q$.
$\tilde K$ together with all these additional points
(which are nonsingular, therefore free of multiplicity 1)
form a cluster $K'$ with Enriques diagram $\D'$,
with $f \in H_{K',\mu'}$. As
$(\D',\mu') \ge (\D,\mu)$, proposition \ref{cideal}
says that there is a cluster $K$ with
Enriques diagram $\D$ such that
$f \in H_{K',\mu'} \subseteq H_{K,\mu}$.
On the other hand
\cite[4.2.7]{Cas00} says that the general
member of $H_{K,\mu}$ defines a germ of type $(\D, \mu)$
so we are done.
%\enlargethispage{\baselineskip}

Let us now prove the only if part, so assume that,
for every $f\in \C[[x,y]]$ defining a reduced germ of
type $(\tilde \D, \tilde \mu)$, there exists an ideal
$I\subset \C[[x,y]]$, with $f \in I$,
whose general member defines a reduced germ of type $(\D, \mu)$.
We first reduce to the case that $I$ has no fixed part.
Indeed, for $n$ big enough and $h \in (x,y)^n$, the
types of $f$ and $f+h$ coincide (see for instance
\cite[7.4.2]{Cas00}), and also the types of $g$ and
$g+h$ for $g$ general in $I$, so we can take $I+ (x,y)^n$
instead of $I$, and this has no fixed part.
Then by \cite[7.2.13]{Cas00} the Enriques diagram of the weighted
cluster $BP(I)$ of base points of $I$ is $(\D, \mu)$ plus some free
vertices of multiplicity one; let $K$ be the subcluster of $BP(I)$
whose Enriques diagram is $\D$. As $f \in I$, $f$ goes through the
weighted cluster $BP(I)$, and therefore $f \in H_{K, \mu}$.
By \cite[4.5.4]{Cas00} this means that
the \emph{value} of $f$ at each point $p \in K$ is at
least $v_p$. Add to the cluster $(\tilde K, \tilde \mu)$
of singular points of $f$ all points on $f=0$ which
belong to $K$, weighted with multiplicity 1
(these are all infinitely near points at which $f=0$
is smooth). Then the resulting cluster $(K',\mu')$
satisfies $H_{K',\mu'} \subseteq H_{K,\mu}$
and by \ref{cideal} we obtain $(\D',\mu') \ge (\D,\mu)$,
where $\D'$ is the Enriques diagram of $K'$.
\end{proof}

%\enlargethispage{6mm}
\begin{figure}
  \begin{center}
    \mbox{\psfig{file=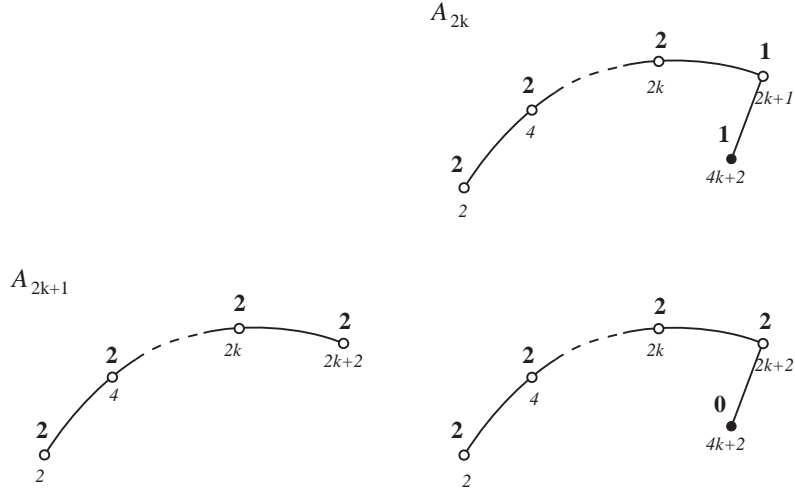}}
%    \mbox{\includegraphics{A.eps}}
    \caption{In white, the vertices of the isomorphic subdiagrams,
      in bold shape, the systems of multiplicities and, in italics,
      the systems of values.
      With notations as in proposition \ref{cideal},
      the top right diagram is $(\D, \mu)$,
      the bottom left is $(\D', \mu')$ and the
      bottom right is $(\D, \nu)$.
      \label{incidA}}
  \end{center}
\end{figure}

\begin{Exa}
Let $A_k$, $D_k$, $E_k$, $J_{k,p}$ and so on denote the
types of germs of curve of Arnold's lists (cf. \cite{Arn76}).
Then for every $k, d>0$, $A_{k+d}$ is linearly adjacent to $A_k$,
$D_{k+d}$ is linearly adjacent to $D_k$,
$E_{k+d}$ is linearly adjacent to $E_k$,
$J_{k+d,p+d}$ is linearly adjacent to $J_{k+d,p}$
and to $J_{k,p+d}$ and so on.
To see this, just take the weighted Enriques
diagrams corresponding to each type, and apply
theorem \ref{linesp}. For instance, figure \ref{incidA}
shows the Enriques diagrams corresponding to types $A_{2k}$
and $A_{2k+1}$ with the corresponding isomorphic subdiagrams,
the multiplicities and the values involved. 
All other cases are handled similarly.
\end{Exa}

%\enlargethispage{12mm}

\begin{Exa}
\label{ex3}
The simplest example in which one needs to consider
$(\D',\mu') \ne (\tilde \D, \tilde \mu)$ is to
prove that a triple point ($D_4$ in Arnold's notation)
is linearly adjacent to the tacnode of type $A_3$ in Arnold's
notation. Indeed, in this case $(\D', \mu')$ is obtained
from the triple point by adding a free point with multiplicity
1 to it (see figure \ref{AD}).
\begin{figure}
  \begin{center}
    \mbox{\psfig{file=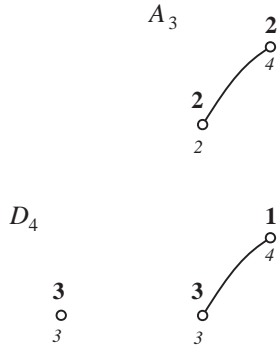}}
%    \mbox{\includegraphics{AD.eps}}
    \caption{Diagrams corresponding to example \ref{ex3}.
          \label{AD}}
  \end{center}
\end{figure}
\end{Exa}

\begin{Rem}
  In the proof of the only if part of \ref{linesp} one just needs to
assume that there exists a $f\in \C[[x,y]]$
defining a reduced germ of type $(\tilde \D, \tilde \mu)$
and an ideal $I\subset \C[[x,y]]$, with $f \in I$,
whose general member defines a reduced germ of type $(\D, \mu)$.
It follows therefore that if a germ of curve of type
$(\tilde \D, \tilde \mu)$ can be deformed linearly to
a germ of type $(\D, \mu)$, then $(\tilde \D, \tilde \mu)$
is linearly adjacent to $(\D, \mu)$.
\end{Rem}

\section{Non-linear adjacency}

%\enlargethispage{\baselineskip}
We have shown in the preceding section a criterion
to decide whether a type
is or is not linearly adjacent to another.
Non-linear adjacencies are a much subtler subject,
as shown by example \ref{contrex} below,
and we cannot give a criterion to decide in all cases.
However, we are able to give a necessary condition
and a sufficient condition.

We say that a weighted Enriques diagram $(\D,\mu)$
is \emph{tame} whenever it is consistent or the
sequence of unloadings that determines, leading to
a consistent Enriques diagram $(\D, \mu')$ is
tame (see \cite[4.7]{Cas00}); the fact that
$(\D,\mu)$ is tame or not depends on the multiplicities
and the proximities between vertices of $\D$, i.e.
on the combinatorial properties of the weighted cluster.
Moreover, $(\D,\mu)$ is \emph{tame} if and only
if for every cluster $K$ with Enriques diagram $\D$ the condition
$\dim \C[[x,y]]/H_{K,\mu}=\deg (\D, \mu)$ holds
(see \cite[4.7.3]{Cas00}).

For every Enriques diagram $\D$, endowed with an
admissible ordering $\preceq$ of the vertices,
there is a variety $\Cl(\D,\preceq)$ parameterizing all
ordered clusters with ordered Enriques diagram $(\D, \preceq)$
(see \cite{Roe?4}). In the sequel we shall
make use of these spaces and the results on
their relative positions in the variety of all clusters
obtained in \cite{Roe?4}. In particular,
we write $(\D, \preceq) \rightsquigarrow (\D', \preceq')$ to mean
$\Cl(\D',\preceq')\subset \overline{\Cl(\D, \preceq)}$.

We begin with a sufficient condition for adjacency.

\begin{Pro}
\label{suf}
Let $(\D,\mu)$, $(\tilde \D,\tilde \mu)$ be types, and
assume that there exist a weighted consistent Enriques diagram
$(\D',\mu')$, differing from $(\tilde \D,\tilde \mu)$ at most
in some free vertices of multiplicity one,
an Enriques diagram $\D_0$ with the same
number of vertices as $\D$, and admissible orderings $\preceq$
and $\preceq_0$ of $\D$ and $\D_0$ respectively satisfying
\begin{enumerate}
\item $(\D,\preceq) \rightsquigarrow (\D_0,\preceq_0)$,
\item $(\D_0,\preceq_0,{\boldsymbol \mu})$ is tame, and
\item $(\D', \mu') \ge (\D_0,\preceq_0,{\boldsymbol \mu})$,
\end{enumerate}
where $\boldsymbol \mu$ is the vector of multiplicities
of $(\D,\mu)$ for the ordering $\preceq$ of $\D$.
Then the type $(\tilde \D,\tilde \mu)$ is
adjacent to the type $(\D,\mu)$.
\end{Pro}
\begin{proof}
 Let $C$ be a germ of curve of type $(\tilde \D,\tilde \mu)$; we
have to see that there is a family of germs containing $C$ whose
general member is of type $(\D,\mu)$. Let $f \in \C[[x,y]]$ be
an equation of $C$, and let $\tilde K$ be the cluster
of singular points of $C$. For each vertex $p$
of $\D'$ not in $\tilde \D$, whose predecessor is denoted by $q$,
choose a point on $C$ on the first neighbourhood
of the point corresponding to the vertex $q$.
$\tilde K$, together with all these additional points
(which are nonsingular, therefore free of multiplicity 1)
form a cluster $K'$ with Enriques diagram $\D'$,
with $f \in H_{K',\mu'}$. As
$(\D',\mu') \ge (\D_0,\preceq_0,{\boldsymbol \mu})$, proposition \ref{cideal}
says that there is a cluster $K_0$ with
Enriques diagram $\D_0$ such that
$f \in H_{K',\mu'} \subseteq H_{K_0,\preceq_0,{\boldsymbol\mu}}$.
The hypothesis $\D \rightsquigarrow \D_0$ says that we can
deform $K_0$ to a family $K_t$ of clusters, $t \in \Delta \subset \C$,
where $\Delta$ is a suitably small disc, such that for $t \ne 0$
the cluster $K_t$ has Enriques diagram $\D$. Now the
$H_{K_t, \mu}$ form a family of linear subspaces of $\C[[x,y]]$
with constant codimension (because $(\D_0,\preceq_0,{\boldsymbol \mu})$ is tame
and $(\D,\mu)$ is consistent) and therefore determine
a family of germs which contain $f$ % \in H_{K_0, \mu}$
and whose general member has type $(\D,\mu)$, as wanted.
\end{proof}

If needed, it is not hard to obtain from the
family described in the proof of proposition \ref{suf}
a one-dimensional family $C_t$ with the desired properties
and $C_0=C$, even explicitly.
For the particular case when $\D$ is unibranched,
the reader may find details on the family
$H_{K_t, \mu}$, with explicit equations,
in \cite[3]{Roe01a}.

Note that, as in the linear case, the interest of proposition
\ref{suf} lies in the fact that the conditions can be
checked directly on the Enriques diagrams, using their
combinatorial properties. This is always true for the conditions
that $(\D_0,\preceq_0,{\boldsymbol \mu})$ is tame, and
$(\tilde \D,\tilde \mu) \ge (\D_0,\preceq_0,{\boldsymbol \mu})$.
The condition $(\D,\preceq) \rightsquigarrow (\D_0,\preceq_0)$
is more difficult to handle, but in some cases (such as when
$\D$ has no satellite points, or when it is unibranched)
it can also be determined from the combinatorial properties of
$\D$ and $\D_0$ (see \cite{Roe?4}) using proximity matrices.

Next we prove a necessary condition
for adjacency (other necessary conditions,
involving invariants such as the codimension
or the Milnor number, are also known).

\begin{Pro}
 Let $(\D,\mu)$, $(\tilde \D,\tilde \mu)$ be types such that
there exists a family of curves $C_t$, $t \in \Delta \subset \C$,
whose general members are of type
$(\D,\mu)$ and with $C_0$ of type $(\tilde \D,\tilde \mu)$.
Then there exist a weighted consistent Enriques diagram
$(\D',\mu')$, differing from $(\tilde \D,\tilde \mu)$ at most
in some free vertices of multiplicity one,
an Enriques diagram $\D_0$ with the same
number of vertices as $\D$ and admissible orderings
$\preceq$ and $\preceq_0$ of $\D$ and $\D_0$
respectively such that
\begin{enumerate}
\item $(\D',\mu') \ge (\D_0,\preceq_0, {\boldsymbol \mu})$,
  and
\item the matrix $P_0^{-1}P$, where $P$ and $P_0$ are the proximity
matrices of $(\D,\preceq)$ and $(\D_0,\preceq_0)$ respectively, has
no negative entries.
\end{enumerate}
\end{Pro}

\begin{proof}
Let $S_t \longrightarrow \Spec \C[[x,y]]$ be a desingularization
of the family $C_t$, $t\ne 0$ (\cite{Zar65II}, see also
\cite{Wah74}). Because of the universal
property of the space $X_{n-1}$ of all ordered clusters of $n$ points (see
\cite{Har85} or \cite{Roe?4}) this induces a family of clusters $K_t$
(parameterized by a possibly smaller punctured disc $\Delta' \setminus \{0\}$)
which can be uniquely extended taking $K_0=\lim_{t \rightarrow 0} K_t$
($X_{n-1}$ is projective and therefore complete). All
clusters of this family except maybe $K_0$ have type $\D$,
and for all $t \in \Delta'$, it is easy to see that
$C_t$ goes through the weighted cluster $(K_t,\mu)$.
Taking $\D_0$ to be the Enriques diagram of $K_0$, both
claims follow (see \cite{Roe?4} for the second claim).
\end{proof}

Obviously this implies
\begin{Cor}
\label{nec}
 Let $(\D,\mu)$, $(\tilde \D,\tilde \mu)$ be types such that
$(\tilde \D,\tilde \mu)$ is adjacent to $(\D,\mu)$.
Then there exist a weighted consistent Enriques diagram
$(\D',\mu')$, differing from $(\tilde \D,\tilde \mu)$ at most
in some free vertices of multiplicity one,
an Enriques diagram $\D_0$ with the same
number of vertices as $\D$, and admissible orderings $\preceq$
and $\preceq_0$ of $\D$ and $\D_0$ respectively such that
\begin{enumerate}
\item $(\D', \mu') \ge (\D_0,\preceq_0,{\boldsymbol\mu})$, and
\item the matrix $P_0^{-1}P$, where $P$ and $P_0$ are the proximity
matrices of $(\D,\preceq)$ and $(\D_0,\preceq_0)$ respectively,
has no negative entries.
\end{enumerate}
\end{Cor}

Again, the interest of
\ref{nec} lies in the fact that the conditions can be
checked directly on the Enriques diagrams, using their
combinatorial properties. Thus we prove, for example, that some
types (including all irreducible curve
singularities with a single characteristic exponent
$m/n$ with $n<m<2n$) allow only linear adjacencies:

\begin{Cor}
   Let $(\D,\mu)$, $(\tilde \D,\tilde \mu)$ be types such that
$(\tilde \D,\tilde \mu)$ is adjacent to $(\D,\mu)$, and suppose that
$\D$ has at most two free vertices. Then $(\tilde \D,\tilde \mu)$ is
linearly adjacent to $(\D,\mu)$.
\end{Cor}
\begin{proof}
  If $p$ is a satellite vertex of $\D$ then there are
at least two vertices in $\D$ preceding it (namely,
the two vertices to which $p$ is proximate). Therefore,
if $\D$ has only one free vertex then it consists of the
root alone, and if it has two free vertices they must
be the root and another vertex which is the unique one
which has the root as immediate predecessor. Under
these conditions, it is not hard to see
that, given any admissible ordering $\preceq$ on $\D$, if
$(\D_0,\preceq_0)$ is an ordered Enriques diagram
such that the matrix $P_0^{-1}P$ has no negative
entries, where $P$ and $P_0$ are the proximity
matrices of $(\D,\preceq)$ and $(\D_0,\preceq_0)$ respectively,
then $(\D,\preceq)=(\D_0,\preceq_0)$. Now the claim
follows from \ref{nec} and \ref{linesp}.
\end{proof}

The fact that the varieties $\Cl(\D)$ do not form a stratification
of the space of all clusters (i.e. there exist $\D$, $\D'$
with $\Cl(\D') \cap \overline {\Cl(\D)} \ne \emptyset$ and
$\D \not\rightsquigarrow \D'$), which is proved in \cite{Roe?4},
implies that the equisingularity classes do not form a stratification
of $\C[[x,y]]$ (i.e. there exist types
$(\D, \mu)$, $(\tilde \D,\tilde \mu)$ and curves of
type $(\tilde \D,\tilde \mu)$ that can be deformed to
curves of type $(\D, \mu)$ without $(\tilde \D,\tilde \mu)$
being adjacent to $(\D, \mu)$). This is shown in the
following example:

\begin{Exa}
  \label{contrex}
Let $(\D, \mu)$, $(\D', \mu')$, $(\tilde \D,\tilde \mu)$ be the
Enriques diagrams of figure \ref{figcontrex}.
In \cite{Roe?4} it is shown that there exist clusters
$K$ and $K'$ with Enriques diagram $\D'$ such that
$K'$ can be deformed to clusters with Enriques diagram
$\D$ and $K$ can not.
If $C$ is a curve of type $(\tilde \D,\tilde \mu)$,
and $(K', \mu')$ is the cluster (of type $(\D' , \mu ')$) formed by the
singular points and the two first nonsingular points on each branch of
$C$, then it is not hard to deform
it to curves of type $(\D, \mu)$, using the method
of the proof of proposition \ref{suf}.
On the other hand, $(\tilde \D,\tilde \mu)$ is
not adjacent to $(\D, \mu)$; this can be proved
using that $K$ cannot be deformed to clusters
with Enriques diagram $\D$ or, more easily, by
observing that both types have the
same codimension.
\begin{figure}
  \begin{center}
    \mbox{\psfig{file=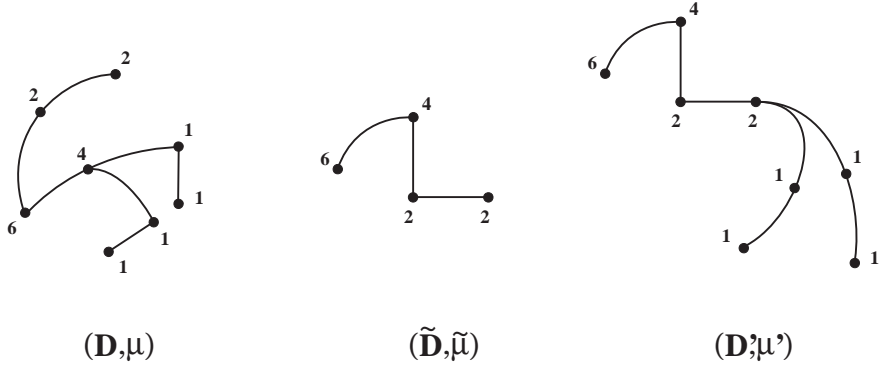}}
%    \mbox{\includegraphics{contrex.eps}}
    \caption{Enriques diagrams corresponding
    to the types of example \ref{contrex}. \label{figcontrex}}
  \end{center}
\end{figure}
\end{Exa}

\section{Non-linear adjacency via Hilbert schemes}
Non-linear adjacency can be approached using Hilbert schemes
instead of varieties of clusters. In fact, it is possible to
give a characterization of all adjacencies in terms of
the relative positions of some subschemes of the Hilbert
scheme of points on a surface. However, these relative
positions are in general not known, so the answer obtained
using Hilbert schemes is theoretical and not easy to put in practice,
in contrast with the criteria given above, which are
combinatorial and can be effectively applied.

As customary, $\Hilb^n R$ will denote the Hilbert
scheme parameterizing ideals of colength $n$ in
$R= \C[[x,y]]$. We consider also the ``nested Hilbert
scheme''
$Z_{n_1, n_2}R \subset (\Hilb^{n_1} R) \times (\Hilb^{n_2} R)$
studied by J. Cheah, which parameterizes pairs of ideals
$(I_1, I_2)$ with $I_1 \supset I_2$
(see \cite{Cheah}, \cite{Che98a}).
For every type $(\D,\mu)$, let $\Hilb_\D^{\mu} R$ be
the subset of $\Hilb^n R$ parameterizing the
ideals $H_{K,\mu}$ where $K$ are clusters
with Enriques diagram $\D$, and $n=\deg (\D, \mu)$. It is known that
$\Hilb_\D^{\mu} R$ is a locally closed irreducible
subscheme of $\Hilb^n R$
(see \cite{KP99}, \cite{NV97}, \cite{Lossen}, for example);
$\overline {\Hilb_\D^{\mu} R}$ will denote its closure
in $\Hilb^n R$.

\begin{Teo}
\label{hilbesp}
Let $(\D,\mu)$, $(\tilde \D,\tilde \mu)$ be types.
$(\tilde \D,\tilde \mu)$ is adjacent to $(\D,\mu)$ if
and only if there exists a weighted consistent Enriques diagram
$(\D',\mu')$, differing from $(\tilde \D,\tilde \mu)$ at most
in some free vertices of multiplicity one, satisfying
$\Hilb^{\mu'}_{\D'} R \subset \pi' \pi^{-1}
\left(\overline {\Hilb^{\mu}_{\D}R}\right)$,
where $\pi$ and $\pi'$ are the projections of
$Z_{n,n'}R$ onto $\Hilb^{n}R$ and $\Hilb^{n'} R$
respectively, and $n=\deg (\D, \mu)$,
$n'=\deg (\D', \mu')$.
\end{Teo}

To prove theorem \ref{hilbesp} we shall use the following
lemma:

\begin{Lem}
\label{unicf}
  Let $(\D,\mu)$, $(\D', \mu')$ be types such that
$(\D', \mu')$ is adjacent to $(\D,\mu)$. Then
for every $f \in \C[[x,y]]$
defining a reduced germ of curve of type $(\D',\mu')$,
there exists an ideal $I \in \overline {\Hilb_\D^{\mu} R}$ with
$f \in I$.
\end{Lem}
\begin{proof}
Let $f \in \C[[x,y]]$ be a germ of equation
of a curve of type $(\D',\mu')$. Because of the adjacency,
there exists a family of germs $f_t$, $t \in \Delta \subset \C$,
whose general members are of type
$(\D,\mu)$ and with $f_0=f$.
Let $S_t \longrightarrow \Spec \C[[x,y]]$ be a desingularization
of the family $f_t$, $t\ne 0$ (\cite{Zar65II}, see also
\cite{Wah74}). Because of the universal
property of the space of all clusters (see \cite{Har85} or
\cite{Roe?4}) this induces a family of clusters $K_t$
(parameterized by a possibly smaller punctured disc
$\Delta' \setminus \{0\}$). Now the
$I_t=H_{K_t, \mu}$ form a (complex) one-dimensional family
inside $\Hilb_\D^{\mu} S$
which can be uniquely extended with
$I_0=\lim_{t \rightarrow 0} I_t$.
It is easy to see that, for all $t \in  \Delta'$,
$f_t \in I_t$, so the claim follows for $I=I_0$.
\end{proof}

\begin{proof}[Proof of theorem \ref{hilbesp}]
The \emph{if} part of the claim is proved in a similar
way to the proof of proposition \ref{suf};
we leave the details to check for the reader.
For the \emph{only if} part of the claim
we shall prove that assuming
$(\tilde \D,\tilde \mu)$ is adjacent to $(\D,\mu)$
and that there exists no consistent weighted Enriques diagram
$(\D',\mu')$, differing from $(\tilde \D,\tilde \mu)$ only
in free vertices of multiplicity one, in the conditions
of the claim, leads to contradiction.

The second assumption means that, for every
consistent Enriques diagram
$(\D',\mu')$, differing from $(\tilde \D,\tilde \mu)$ only
in free vertices of multiplicity one, there are clusters $K'$
with $H_{K',\mu'} \in \Hilb^{\mu'}_{\D'} R \setminus
\pi' \pi^{-1} \left(\overline {\Hilb^{\mu}_{\D}R}\right)$.
Consider the sequence of weighted Enriques diagrams
defined as follows. $(\D_1, \mu_1)$ is obtained
from $(\tilde \D, \tilde \mu)$ by adding
$$\tilde \mu_p- \sum_{q \text{ prox. to } p} \tilde \mu_q$$
free successors of multiplicity 1 to each $p \in \tilde \D$,
and for $k>1$, $(\D_k, \mu_k)$ is obtained from
$(\D_{k-1}, \mu_{k-1})$ by adding a free successor
of multiplicity 1 to
each extremal vertex (which will be free of multiplicity 1).
Obviously $(\D_{k-1}, \mu_{k-1})$ is a subdiagram
of $(\D_k, \mu_k)$ for all $k>1$, and
it is not hard to see that the
map $F_k:\Hilb^{\mu_k}_{\D_k}R \longrightarrow
\Hilb^{\mu_{k-1}}_{\D_k-1}R$
defined by sending $H_{K,\mu_k}$ to $H_{\breve K, \mu_{k-1}}$,
where $\breve K$ is the subcluster of $K$ with diagram
$\D_{k-1}$, satisfies
$$F_k\left(\Hilb^{\mu_k}_{\D_k} R \setminus
\pi_k \pi^{-1} \left(\overline {\Hilb^{\mu}_{\D}R}\right)\right)=
\Hilb^{\mu_{k-1}}_{\D_{k-1}} R \setminus
\pi_{k-1}\pi^{-1} \left(\overline {\Hilb^{\mu}_{\D}R}\right).$$
Therefore we can construct a sequence of clusters
$K_1, K_2, \dots$ such that
$H_{K_k, \mu_k} \in \Hilb^{\mu_k}_{\D_k}R \setminus
\pi_k \pi^{-1} \left(\overline {\Hilb^{\mu}_{\D}R}\right)$
and each $K_k$ is obtained from
$K_{k-1}$ by adding in the first neighbourhood of each extremal point
a free point of multiplicity one. But then there exists
a reduced germ $f$ of type $(\tilde \D, \tilde \mu)$
belonging to all $H_{K_k, \mu_k}$ (see \cite[5.7]{Cas00}]).
%and in fact it is easy to see that $\cap H_{K_k, \mu_k}=(f)$.
Now by lemma \ref{unicf}, there exists an ideal
$I \in \overline {\Hilb_\D^{\mu} R}$ with
$f \in I$; as $\dim_{\C} \C[[x,y]]/I = n$,
we must have $I \supset (x,y)^n$ also.
On the other hand, applying \cite[5.7.1]{Cas00}
and \cite[7.2.16]{Cas00},
for $k$ big enough we infer that
$H_{K_k, \mu_k} \subset (f)+(x,y)^n$,
which implies $H_{K_k, \mu_k} \subset I$, a contradiction.
\end{proof}

\begin{Rem}
Linear adjacencies may also be dealt with using Hilbert
schemes; indeed, with notations as above,
$(\tilde \D,\tilde \mu)$ is linearly adjacent to $(\D,\mu)$ if
and only if there exists a weighted consistent Enriques diagram
$(\D',\mu')$, differing from $(\tilde \D,\tilde \mu)$ at most
in some free vertices of multiplicity one, satisfying
$\Hilb^{\mu'}_{\D'} R \subset \pi' \pi^{-1}
\Hilb^{\mu}_{\D}R$. Again this criterion is hard to
apply, in contrast to the purely combinatorial we gave
before. We skip the proof, which adds no new ideas
to what we did before.
\end{Rem}

\begin{Rem}
For types $(\D,\mu)$ where $\D$ has three vertices or
less, the closure of $\Hilb^\mu_\D R$ is known, due to
  the works \cite{Evain} and \cite{eva?1}  of \'Evain; so in this case
the Hilbert scheme method does give a characterization of adjacencies.
Very few other particular situations can be handled
explicitly; we would like to mention an example due
to Russell (see \cite{Rus??}) in which the study of
the Hilbert scheme provides an example (like \ref{contrex})
showing that types do not stratify $\C[[x,y]]$.
\end{Rem}

\bibliographystyle{amsplain}
\bibliography{Biblio}

\providecommand{\bysame}{\leavevmode\hbox to3em{\hrulefill}\thinspace}
\begin{thebibliography}{10}

\bibitem{Arn76}
V.~I. Arnold, \emph{Local normal forms of functions}, Invent. Math. \textbf{35}
  (1976), 87--109.

\bibitem{Cas00}
E.~Casas-Alvero, \emph{Singularities of plane curves}, London Math. Soc.
  Lecture Notes Series, no. 276, Cambridge University Press, 2000.

\bibitem{Cheah}
J.~Cheah, \emph{The cohomology of smooth nested {H}ilbert schemes of points},
  Ph.D. thesis, Chicago, 1994.

\bibitem{Che98a}
\bysame, \emph{The virtual hodge polynomials of nested {H}ilbert schemes and
  related varieties}, Math. Z. \textbf{227} (1998), no.~3, 479--504.

\bibitem{EC15}
F.~Enriques and O.~Chisini, \emph{Lezioni sulla teoria geometrica delle
  equazioni e delle funzioni algebriche}, N. Zanichelli, Bologna, 1915.

\bibitem{Evain}
L.~{\'E}vain, \emph{Collisions de trois gros points sur une surface
  alg{\'e}brique}, Ph.D. thesis, Universit{\'e} de {N}ice, 1997.

\bibitem{eva?1}
\bysame, \emph{Compactification of configuration spaces via {H}ilbert schemes},
  preprint (2001), \\http://xxx.lanl.gov/abs/math/0107041.

\bibitem{GSG92}
A.~Granja and T.~S\'anchez-Giralda, \emph{Enriques graphs of plane curves.},
  Comm. Algebra \textbf{20} (1992), no.~2, 527--562.

\bibitem{Har85}
B.~Harbourne, \emph{Complete linear systems on rational surfaces}, Trans.
  A.M.S. \textbf{289} (1985), 213--226.

\bibitem{KP99}
S.~Kleiman and R.~Piene, \emph{Enumerating singular curves on surfaces}, Proc.
  Conference on {A}lgebraic {G}eometry: {H}irzebruch 70 (Warsaw 1998), vol.
  241, A.M.S. Contemp. Math., 1999, pp.~209--238.

\bibitem{Lossen}
C.~Lossen, \emph{The geometry of equisingular and equianalytic families of
  curves on a surface}, Ph.D. thesis, Universit\"at Kaiserslautern, 1998.

\bibitem{NV97}
A.~Nobile and O.~Villamayor, \emph{Equisingular stratifications associated to
  families of planar ideals}, J. Alg. \textbf{193} (1997), 239--259.

\bibitem{Roe?4}
J.~Ro\'e, \emph{Varieties of clusters and {E}nriques diagrams}, preprint,
  \\http://xxx.lanl.gov/abs/math/0108023.

\bibitem{Roe01a}
\bysame, \emph{On the conditions imposed by tacnodes and cusps}, Trans. A.M.S.
  \textbf{353} (2001), no.~12, 4925--4948.

\bibitem{Rus??}
H.~Russell, \emph{Counting singular plane curves via {H}ilbert schemes},
  Preprint (2000), \\http://xxx.lanl.gov/abs/math/0011214.

\bibitem{Sie77}
D.~Siersma, \emph{Periodicities in {A}rnold's lists of singularities}, Real and
  Complex Singularities, {O}slo 1976 (P.~Holm, ed.), Sijthoff \& Noordhoof,
  1977, pp.~497--524.

\bibitem{Wah74}
J.~Wahl, \emph{Equisingular deformations of plane algebroid curves}, Trans.
  A.M.S. \textbf{193} (1974), 143--170.

\bibitem{Zar32}
O.~Zariski, \emph{On the topology of algebroid singularities}, Amer. J. Math.
  \textbf{54} (1932), 433--465.

\bibitem{Zar65I}
\bysame, \emph{Studies in equisingularity {I}}, Amer. J. Math. \textbf{87}
  (1965), 507--536.

\bibitem{Zar65II}
\bysame, \emph{Studies in equisingularity {II}}, Amer. J. Math. \textbf{87}
  (1965), 972--1006.

\end{thebibliography}

\end{document}